\documentclass[12pt]{article}
\usepackage{amsthm,amsmath,amssymb,amsfonts,latexsym}

\newtheorem{theorem}{\qquad Theorem}
\newtheorem{proposition}{\qquad Proposition}
\newtheorem{corollary}{\qquad Corollary}

\title{Expansions for solutions of the Schlesinger equation at a singular point\footnote{The work is supported by grants RFBR 12-01-33058-mol\_a\_ved, RFBR 12-01-31414-mol\_a and Simons-IUM fellowship-2012}}

\author{I.\,V. Vyugin \\ {\small Institute for Information Transmission Problems RAS} \\ {\small vyugin@gmail.com}}

\date{}

\begin{document}
\maketitle

\begin{abstract}
A local behavior of solutions of the Schlesinger equation is studied. We obtain expansions for this solutions, which converge in some
neighborhood of a singular point. As a corollary the similar result for the sixth Painlev\'e equation was obtained. In our analysis, we use the
isomonodromic approach to solve this problem.
\end{abstract}

Keywords: Isomonodromic deformation, Fuchsian system, Schlesinger equation, bundle with connection

MSC 34M56, 34M55, 34M03

\section{Introduction}

We study a local behavior of solutions of the Schlesinger equation. We present solutions of this equation in the form of power series or
logarithmic-power series. This series are converge in some neighborhood of a singular point. As a corollary we obtain a similar result for
description of the behavior of solutions of the sixth Painlev\'e equation in some sectorial neighborhood. We use the isomonodromic approach to
solve this problem.

Let us consider the following system of analytical partial differential equations
\begin{eqnarray}\label{Schlesinger}
dB_i=-\sum_{j=1,j\not=i}^n\frac{[B_i,B_j]}{a_i-a_j}d(a_i-a_j),\qquad i=1,\ldots,n,
\end{eqnarray}
where $B_i$ ($i=1,\ldots,n$) --- are analytical $p\times p$-matrix functions of the variable $a=(a_1,\ldots,a_n)$, $[B_i,B_j]$ denotes the
commutator of matrices $B_i$ and $B_j$. The matrix-functions $B_i(a)$ are defined and meromorphic (see B. Malgrange \cite{Mal}, R. Gontsov and
I. Vyugin \cite{VG2}) on the space
$$
\{a|a=(a_1,\ldots,a_n)\in \overline{\mathbb{C}}^n\setminus \bigcup_{i,j}A_{ij}\},\qquad A_{ij}=\{a\,|\, a_i= a_j\}.
$$
This system is called {\it Schlesinger equation} (read more in A.A. Bolibruch \cite{Bo2}). Divisor of the Schlesinger equation is the following
set
\begin{eqnarray*}
\Omega=\bigcup_{i,j}A_{ij}.
\end{eqnarray*}
We are going to describe a local form of solutions of Schlesinger equation (\ref{Schlesinger}) in a neighborhood of the point
$a^0=(a_1^0,\ldots,a_n^0)$, which belongs to the following singular set
\begin{eqnarray*}
a^0\in\Omega'=\Omega\setminus\left(\bigcup_{i,j,k}A_{ijk}\right),\qquad A_{ijk}=\{a\,|\, a_i= a_j=a_k\}.
\end{eqnarray*}
We obtain the local expansions of the solutions of the system (\ref{Schlesinger}) in the form of power and logarithmic-power series of
$(a_s-a_r)$ (if $a^0\in A_{sr}$), which converges in some neighborhood of the point $z=a^0$ (the first version of these results see \cite{V}).
These series have terms of complex degrees.
\begin{theorem}\label{th1}
Any solution of two dimensional Schlesinger equation (\ref{Schlesinger}) can be represented in the neighborhood of a point
$a^0=(a_1^0,\ldots,a_n^0)\in\Omega'$, where $a_r^0=a_s^0$, $r\not=s$, in one of two following forms:
\begin{itemize}
\item{}
$b_{kl}^i(a)=F_1^{kli}(a)+(a_s-a_r)^{\varphi}F_2^{kli}(a)+(a_s-a_r)^{-\varphi}F_3^{kli}(a)$, $\varphi\in\mathbb{C}$ in the general case;
\item{}
$ b_{kl}^i(a)=F_1^{kli}(a)+F_2^{kli}(a)\ln (a_s-a_r)+F_3^{kli}(a)\ln^2 (a_s-a_r)$ in the degenerate case,
\end{itemize}
where $F_1^{kli}(a),F_2^{kli}(a),F_3^{kli}(a)$ are meromorphic (holomorphic in the generic case) functions, $i=1,\ldots,n$, and $k,l\in\{1,2\}$.
\end{theorem}

The notions of ``general case'' and ``non-general case'' are explained below. Notice that the measure of the systems of non-general case is
equal to zero.

Now consider the case $n=4$, $p=2$, which is equivalent to case of the sixth Painlev\'e equation (\ref{P6}). Without loss of generality, let us
fix three variable $a_1=0$, $a_2=1$, $a_3=\infty$ and denote $a_4$ by $t$. We obtain the system of ordinary differential equations with variable
$t$ and unknown matrix-functions
$$
B_i(t)=\left(\begin{array}{cc}
b_{11}^i(t) & b_{12}^i(t) \\
b_{21}^i(t) & b_{22}^i(t)
\end{array}\right),\qquad i=0,t,1,\infty.
$$
With restrictions above the following corollary holds.

\begin{corollary}\label{Cor}
Any solution of the Schlesinger equation under the above constraints can be represented in the neighborhood of $t=0$ in one of two forms:
\begin{itemize}
\item{}
$b_{kl}^i(t)=F_1^{kli}(t)+t^{\varphi}F_2^{kli}(t)+t^{-\varphi}F_3^{kli}(t)$, $\varphi\in\mathbb{C}$ in the general case;
\item{}
$b_{kl}^i(t)=F_1^{kli}(t)+F_2^{kli}(t)\ln t+F_3^{kli}(t)\ln^2 t$ in the degenerate case,
\end{itemize}
where $F_1^{kli}(t),F_2^{kli}(t),F_3^{kli}(t)$ are meromorphic in $t=0$ functions, $i=0,t,1,\infty$, and $k,l\in\{1,2\}$.
\end{corollary}

Note that the well-known sixth Painlev\'e equation
\begin{eqnarray}\label{P6}
\frac{d^2 w}{dt^2}=\frac{1}{2}\left(\frac{1}{w}+\frac{1}{w-1}+\frac{1}{w-t}\right)\left(\frac{d
w}{dt}\right)^2-\left(\frac{1}{t}+\frac{1}{t-1}+\frac{1}{w-t}\right)\frac{d w}{dt}+\\
+\frac{w(w-1)(w-t)}{t^2(t-1)^2}\left(\alpha+\beta\frac{t}{w^2}+\gamma\frac{t-1}{(w-1)^2}+\delta\frac{t(t-1)}{(w-t)^2}\right),\nonumber
\end{eqnarray}
$\alpha,\,\beta,\,\gamma,\,\delta\in\mathbb{C}$ is equivalent to the system (\ref{Schlesinger}), where
\begin{eqnarray}\label{P6-sol-iso}
w(t)=\frac{tb_{12}^0}{(t+1)b_{12}^0+tb_{12}^1+b_{12}^t}.
\end{eqnarray}
Corollary \ref{Cor} and (\ref{P6-sol-iso}) give the power expansions for solutions of the sixth Painlev\'e equation. A different asymptotics for
sixth Painlev\'e equation was obtained in D. Guzzetti \cite{Guz}, A. Bruno and I. Goryuchkina \cite{BG}, M. Mazzocco \cite{Maz1} and others.

For the sixth Painlev\'e equation, we have an analogue of Corollary \ref{Cor}.

\begin{corollary}\label{P6th}
Any solution $w(t)$ of sixth Painlev\'e equation (\ref{P6}) in the intersection of the given sector for $t$ sufficiently close to singular point
$t=0,1,\infty$ can be represented as a converged power series or as a converged logarithmic-power series:
\begin{itemize}
\item{}
if $G_1G_{\infty}$ is digonalizable, then $w(t)=S(t,t^{\varphi},t^{-\varphi})$, where
$\varphi=\varphi(\alpha,\beta,\gamma,\delta,t_0,w(t_0),w'(t_0))$ can be found approximately;
\item{}
if $G_1G_{\infty}$ is a Jordan block, then $w(t)=S(t,\ln t,\ln^{-1}t)$.
\end{itemize}
\end{corollary}

Using the expressions for $b_{12}^i(t)$ we obtain the following expressions for $w(t)$:
$$
w(t)=\frac{f_1(t)+t^{\varphi}f_2(t)+t^{-\varphi}f_3(t)}{g_1(t)+t^{\varphi}g_2(t)+t^{-\varphi}g_3(t)}
$$
and
$$
w(t)=\frac{f_1(t)+f_2(t)\ln t+f_3(t)\ln^2 t}{g_1(t)+g_2(t)\ln t+g_3(t)\ln^2 t},
$$
where $f_j$, $g_j$, $j=1,2,3$ are meromorphic functions. The denominators of these ratios are not identically zero. We can express these ratios
as power series with powers of $t,$ $t^{\varphi}$, $t^{-\varphi}$, $\ln t$. These power series will be converge in sectorial neighborhood with
any angle $\psi$ and with radius $r$, which depends of $\psi$, $r=r(\psi)$. This sectorial neighborhood is described by the condition: the
denominator of the ratio does not vanishes.

\section{Schlesinger equation and isomonodromic deformations}

In this section we give a description of the Schlesinger equation (\ref{Schlesinger}) as an isomonodromy condition for a family of Fuchsian
systems. Let us consider a Fuchsian system
\begin{eqnarray}\label{Fuks-syst}
\frac{dy}{dz}=\left(\sum_{i=1}^n\frac{B_i^0}{z-a_i^0}\right)y, \qquad B_i^0\in{\rm Mat}_{p\times p}(\mathbb{C}),\quad y(z)\in\mathbb{C}^p.
\end{eqnarray}
The family of such systems
\begin{eqnarray}\label{Fuks-syst-sem}
\frac{dy}{dz}=\left(\sum_{i=1}^n\frac{B_i(a)}{z-a_i}\right)y
\end{eqnarray}
is called {\it isomonodromic} if the following conditions hold:
\begin{itemize}
\item{}
$B_i(a)$ are continuous matrix-functions of $a=(a_1,\ldots,a_n)$;
\item{}
The Fuchsian system (\ref{Fuks-syst-sem}) with any fixed $a$ has fixed monodromy representation
\begin{eqnarray}\label{repr}
\chi: \pi_1(\overline{\mathbb{C}}\setminus\{a_1,\ldots,a_n\},z_0)\longrightarrow GL(p,\mathbb{C}).
\end{eqnarray}
\end{itemize}

{\it Schlesinger isomonodromic family} is a family defined by the equation (\ref{Schlesinger}). An isomonodromic fundamental matrix $Y(z,a)$ of
the Schlesinger isomonodromic family (\ref{Fuks-syst-sem}) satisfies the following condition
\begin{eqnarray*}
Y(\infty,a)\equiv Y(\infty,a_0).
\end{eqnarray*}
The initial data of such family are the coefficients $B_i(a^0)=B_i^0$, $i=1,\ldots,n$ of system (\ref{Fuks-syst}). It is known that the
solutions of Schlesinger equation are meromorphic functions on the space $a\in\overline{\mathbb{C}}^n\setminus\Omega$.

\bigskip

Let us consider the Painlev\'e VI case ($n=4,p=2$, $a_1=0$, $a_2=1$, $a_3=\infty$, $a_4=t$). Usually the following family
\begin{eqnarray}\label{syst-fuks-4}
\frac{dy}{dz}=\left(\frac{B_0(t)}{z}+\frac{B_t(t)}{z-t}+\frac{B_1(t)}{z-1}\right)y
\end{eqnarray}
is considered, where
$$
{\rm tr}B_0={\rm tr}B_t={\rm tr}B_1={\rm tr}B_{\infty}=0,\qquad B_{\infty}=-(B_0+B_t+B_1),
$$
and the matrices $B_0$, $B_t$, $B_1$, $B_{\infty}={\rm diag}(\delta,-\delta)$ are diagonalizable.

The formula (\ref{P6-sol-iso}) gives a solution of sixth Painlev\'e equation (\ref{P6}) with the following constants
\begin{eqnarray*}
\alpha=\frac{(2\lambda_{\infty}-1)^2}{2},\quad\beta=-2\lambda_0^2,\quad\gamma=2\lambda_1^2,\quad\delta=\frac{1}{2}-2\lambda_t^2,
\end{eqnarray*}
where $\lambda_0$, $\lambda_t$,$\lambda_1$, $\lambda_{\infty}$ are eigenvalues of matrices $B_0$, $B_t$, $B_1$, $B_{\infty}$.

\section{Sketch of the proof}

At first, we study the simplest isomonodromic family. Consider the following family of Fuchsian systems
\begin{eqnarray}\label{Norm-sem1}
\frac{dy}{dz}=\left(\frac{B_0}{z}+\frac{B_t}{z-t}\right)y,\qquad B_0,B_t\equiv {\rm const}.
\end{eqnarray}
It is easy to see that it is an isomonodromic family. The systems of this family are mutually equivalent by a linear mapping of $z$. This family
is non-Schlesinger if $B_{\infty}=-B_0-B_t\not= 0$, but this family can be transformed to a Schlesinger family by the gauge transformation
$\tilde{y}=t^{B_{\infty}}y$. The transformed family has the following form
\begin{eqnarray}\label{Norm-sem2}
\frac{dy}{dz}=\left(\frac{B_0'(t)}{z}+\frac{B_t'(t)}{z-t}\right)y,\qquad B_i'=t^{B_{\infty}}B_it^{-B_{\infty}},\quad i=0,t.
\end{eqnarray}
We call this family {\it canonical normalized family}.

Let us write the coefficients $B_0',B_t'$ explicitly. There are two cases: in the first case when the matrix $B_{\infty}$ is diagonalizable, and
in the second case $B_{\infty}$ is a Jordan block:
\begin{itemize}
\item{}
First case, $B_{\infty}=\left(\begin{array}{cc}
\lambda_1 & 0    \\
0 & \lambda_2
\end{array}\right)$ and
$ B_i'=\left(\begin{array}{cc}
b_{11}^i & b_{12}^i t^{\lambda_1-\lambda_2} \\
b_{21}^i t^{\lambda_2-\lambda_1} & b_{22}^i
\end{array}\right).
$
\item{}
Second case, $B_{\infty}=\left(\begin{array}{cc}
\lambda & 1    \\
0 & \lambda
\end{array}\right)$ and
$ B_i'=\left(\begin{array}{cc}
b_{11}^i+b_{21}^i \ln t & b_{12}^i+(b_{22}^i-b_{11}^i) \ln t-b_{21}^i \ln^2 t \\
b_{21}^i & b_{22}^i-b_{21}^i \ln t.
\end{array}\right).
$
\end{itemize}

\bigskip

Now we study the limit of the family (\ref{Norm-sem2}) as $t\to 0$. We would like to find a limit Fuchsian system. For the existence of this
limit, we impose the following condition on the real part of $\lambda_1-\lambda_2$
\begin{eqnarray}\label{condB}
\left|{\rm Re}(\lambda_1-\lambda_2)\right|<1.
\end{eqnarray}
The condition (\ref{condB}) implies the following equation (see \cite{Bo2})
$$
\frac{B_0'(t)}{z}+\frac{B_t'(t)}{z-t}=t^{B_{\infty}}\left(\frac{B_0}{z}+\frac{B_t}{z-t}\right)t^{-B_{\infty}}= t^{B_{\infty}}\left(\frac{
B_0+B_t}{z}+O(t)\right)t^{-B_{\infty}}=\frac{B_{\infty}}{z}+o(1).
$$
The limit system as $t\to 0$ of the family (\ref{Norm-sem2}) under the condition (\ref{condB}) is
\begin{eqnarray}\label{Fuks-eul}
\frac{dy}{dz}=\frac{B_{\infty}}{z}y.
\end{eqnarray}
The proof is similar to that given in A.A. Bolibruch \cite{Bo2}.

By the Riemann--Hilbert theory gives that for almost all monodromy representations (\ref{repr}) with generators $G_0'$, $G_1'$, $G_{\infty}'$
there exists a Fuchsian system (\ref{Norm-sem1}) with this monodromy data and given asymptotics (see \cite{AB}). In all other cases, we can
construct such system with one regular singular point (see \cite{VG}). This cases were called in Theorem \ref{th1} ``general case'' and
``non-general case''.

Now let us consider the non-general case. It is the case, when there isn't a system (\ref{Norm-sem1}) having the given monodromy (\ref{repr}),
and exponents in the points $z=0,t,\infty$. The results of I. Vyugin and R. Gontsov \cite{VG} states that there exist the regular system
\begin{eqnarray*}
\frac{dy}{dz}=\left(\frac{B_{-r}}{z^{r+1}}+\ldots+\frac{B_0}{z}+\frac{B_{t}}{z-t_0}\right)y,
\end{eqnarray*}
having the given monodromy (\ref{repr}), and exponents $\beta_t^1$, $\beta_t^2$, $\beta_{\infty}^1$, $\beta_{\infty}^2$ in points $z=t,\infty$,
and $r<3\max(\beta_t^1-\beta_t^2,\beta_{\infty}^1-\beta_{\infty}^2)$. Note that the family
\begin{eqnarray}\label{Syst-reg-sem}
\frac{dy}{dz}=\left(t^r\frac{B_{-r}'}{z^{r+1}}+t^{r-1}\frac{B_{-(r-1)}'}{z^{r}}+\ldots+\frac{B_0'}{z}+\frac{B_t'}{z-t}\right)y,
\end{eqnarray}
where $B_i'=t^{B_{\infty}}B_it^{-B_{\infty}}$, $B_{\infty}=-(B_0+B_t)$, is an isomonodromic. The limit (\ref{Syst-reg-sem}) at $t\to 0$ is
(\ref{Fuks-eul}).

\bigskip

We will use the family (\ref{Norm-sem2}) for the proof of the theorem \ref{th1}.

Let us consider a family of holomorphic vector bundles with logarithmic connection having the following description
\begin{eqnarray*}
(F_t,\nabla_t)=(D_0,D_{\infty},g_{0\infty}^t(z),\omega_0^t,\omega_{\infty}),
\end{eqnarray*}
where $D_0,D_{\infty}$ --- circles with centers $0$ and $\infty$, which has a nonempty intersection $K=D_0\cap D_{\infty}$, $g_{0\infty}^t(z)$
is a holomorphic cocycle $g_{0\infty}^t(z):K\to {\rm GL}(p,\mathbb{C})$ and $\omega_0^t$, $\omega_{\infty}$ is a differential 1-forms of
logarithmic connection $\nabla_t$.

Define the pairs $(F_t,\nabla_t)$ by the following description:
\begin{itemize}
\item{}
$\omega_{\infty}$ is a 1-form of coefficients of the initial system (\ref{syst-fuks-4}), when $t=t_0$, which has the monodromy representation
(\ref{repr}) and generators $G_0,G_t,G_1,G_{\infty}$;
\item{}
$\omega_0^t$ is a 1-form of coefficients of system (\ref{Norm-sem2}), with monodromy (\ref{repr}) with generators $G_0',G_t',G_{\infty}'$
($G_0'=G_0$, $G_t'=G_t$, $G_{\infty}'=G_1G_{\infty}$);
\item{}
cocycle $g_{0\infty}^t(z)$ is a ratio $g_{0\infty}^t(z)=Y_0^t(z)Y_{\infty}^{-1}(z)$, where $Y_0^t(z)$ and $Y_{\infty}(z)$ are fundamental
matrices of the systems
$$
dy=\omega_0^t y,\qquad dy=\omega_{\infty} y,
$$
normalized in $z=\infty$.
\end{itemize}
\begin{proposition}\label{Bun-lim}
Assume that the family
\begin{eqnarray*}
(F_t,\nabla_t)=(D_0,D_{\infty},g_{0\infty}^t(z),\omega_0^t,\omega_{\infty}^t),
\end{eqnarray*}
holomorphically depends on the valuable $t$. If the limit
\begin{eqnarray*}
(F_0,\nabla_0)=(D_0,D_{\infty},g_{0\infty}^0(z),\omega_0^0,\omega_{\infty}^0),\quad t\to 0
\end{eqnarray*}
exists and it is a trivial bundle with trivialization $V^0(z)$, $W^0(z)$, then bundles $F_t$ are trivial bundles, for sufficiently small $t$,
and their trivializations $V^t(z)$, $W^t(z)$ have the limits
\begin{eqnarray*}
\lim_{t\to 0}V^t(z)=V^0(z),\qquad \lim_{t\to 0}W^t(z)=W^0(z),
\end{eqnarray*}
which are uniform for $z\in D_0\cap D_{\infty}$
\end{proposition}

The family (\ref{syst-fuks-4}) can be represented in the neighborhood $D_0$ as
\begin{eqnarray}\label{koeff-0}
\frac{dy}{dz}=\left(V(t,z)\omega_0^t V^{-1}(t,z)+\frac{dV}{dz}V^{-1}\right)y,\qquad V(t,z)=V_{0}(t)+V_1(t)z+\ldots,
\end{eqnarray}
and in the neighborhood $D_{\infty}$ as
\begin{eqnarray}\label{koeff-inf}
\frac{dy}{dz}=\left(W(t,z)\omega_{\infty}W^{-1}(t,z)+\frac{dW}{dz}W^{-1}\right)y,\qquad W(t,z)=W_{0}(t)+W_1(t)z+\ldots.
\end{eqnarray}
We obtain that in the general case the matrix-coefficients of the system (\ref{syst-fuks-4}) are
$$
B_i=V_0(t)B_i'(t)V_0^{-1}(t),\qquad i=0,t,
$$
and
$$
B_i=W_0(t)B_i^0(t)W_0^{-1}(t),\qquad i=1,\infty.
$$
From this formulas and from the form (\ref{Norm-sem2}) we obtain Theorem \ref{th1} and Corollary \ref{Cor} in the general case.

\end{document}